\documentclass{gtart}    

\def\ifplaintex{\expandafter\ifx\csname documentclass\endcsname\relax}

\ifplaintex 
\else
\fi

\def\gtm{{\mathsurround=0pt\it $\cal G\mskip-2mu$eometry \&\ 
$\cal T\!\!$opology $\cal M\mskip-1mu$onographs}}    %  for monographs

\def\gtp{{\mathsurround=0pt\it $\cal G\mskip-2mu$eometry \&\ 
$\cal T\!\!$opology $\cal P\!$ublications}}  % GT publications

\def\recd{{\small Received:\qua\receiveddate\ifx\reviseddate\relax
\else\qquad Revised:\qua\reviseddate\fi\par}} 

\def\volumenumber#1{\def\thevolumenumber{#1}}
\def\volumeyear#1{\def\thevolumeyear{#1}}
\def\volumename#1{\def\thevolumename{#1}}
\def\papernumber#1{\def\thepapernumber{#1}}
\def\pagenumbers#1#2{\def\startpage{#1}\def\finishpage{#2}}
\def\published#1{\def\publishdate{#1}}
\def\received#1{\def\receiveddate{#1}}
\def\revised#1{\def\reviseddate{#1}}
\def\accepted#1{\def\accepteddate{#1}}

\long\def\asciiabstract#1{\long\def\theasciiabstract{#1}}

\let\\\par
\let\thevolumenumber\relax\let\thepapernumber\relax
\let\thevolumeyear\relax\let\startpage\relax
\let\finishpage\relax\let\publishdate\relax\let\receiveddate\relax
\let\reviseddate\relax\let\accepteddate\relax\let\theasciititle\relax
\let\theasciiauthors\relax
\let\theasciiabstract\relax

\let\theerratum\relax\let\theasciiemail\relax
\let\theshortauthors\relax\let\theshorttitle\relax

\def\startpage{1}\def\finishpage{15}\def\thepapernumber{77}

\volumenumber{2}
\volumename{Proceedings of the Kirbyfest}
\volumeyear{1999}

\long\def\maketitlep{   % start of definition of \maketitlep

\count0=\startpage

\gtm\nl        %   GT mongraphs (top left) 
{\small Volume \thevolumenumber: \thevolumename\nl 
\ifx\theerratum\relax\else Erratum \erratumnumber\nl\fi
Pages \startpage--\finishpage\nl}

\vglue 0.1truein   % top margin

{\parskip=0pt\leftskip 0pt plus 1fil\def\\{\par\smallskip}{\ifplaintex\large
\else\Large\fi\bf\thetitle}\par\medskip}   
\vglue 0.05truein 

{\parskip=0pt\leftskip 0pt plus 1fil\def\\{\par}{\sc\theauthors}
\par\medskip}%
 
\vglue 0.03truein 

{\small\leftskip 25pt\rightskip 25pt{\bf Abstract}\stdspace\theabstract

{\bf AMS Classification}\stdspace\theprimaryclass
\ifx\thesecondaryclass\relax\else; \thesecondaryclass\fi\par
{\bf Keywords}\stdspace \thekeywords\par}\vglue 7pt

}   % end of definition of \maketitlep

\font\phead=cmsl9 scaled 950
\font=cmsl9 scaled 1050
\font\pnum=cmbx10 scaled 913
\font\lnum=cmbx10 
\font\pfoot=cmsl9 scaled 950
\font=cmsl9 scaled 1050
\ifplaintex
\headline{\vbox to 0pt{\vskip -4.5mm\line{\small\phead\ifnum
\count0=\startpage ISSN 1464-8997 (on line)
1464-8989 (printed) \hfill {\pnum\folio}\else\ifodd\count0\def\\{ }% 
\ifx\theshorttitle\relax\thetitle\else\theshorttitle\fi\hfill{\pnum\folio}
\else\def\\{ and }{\pnum\folio}\hfill\ifx\theshortauthors\relax\theauthors
\else\theshortauthors\fi\fi\fi}\vss}}
\footline{\vbox to 0pt{\vglue 0mm\line{\small\pfoot\ifnum\count0=\startpage
Published \publishdate:\qua\copyright\ \gtp\hfill\else
\gtm, Volume \thevolumenumber\ (\thevolumeyear)\hfill\fi}\vss
}}
\else
\makeatletter
\def\@oddhead{{\small\lhead\ifnum\count0=\startpage ISSN 1464-8997 (on line)
1464-8989 (printed) \hfill {\lnum\number\count0}\else\ifodd\count0
\def\\{ }\ifx\theshorttitle\relax \thetitle \else\theshorttitle\fi\hfill
{\lnum\number\count0}\else\def\\{ and }{\lnum\number\count0}
\hfill\ifx\theshortauthors\relax 
\theauthors\else\theshortauthors\fi\fi\fi}}\def\@evenhead{@oddhead}
\def\@oddfoot{\small\lfoot\ifnum\count0=\startpage Published \publishdate:\qua\copyright\ \gtp\hfill\else
\gtm, Volume \thevolumenumber\ (\thevolumeyear)\hfill\fi}
\def\@evenfoot{@oddfoot}
\makeatother
\fi

\let\maketitlepage\maketitlep
\let\makeshorttitle\maketitlepage
\let\maketitle\maketitlepage

\newwrite\gtoutfile
\long\gdef\makeheadfile{  %%% start of definition of \makeheadfile
{\def\\{, }\def\s{ }
\immediate\openout\gtoutfile head.xxx
\immediate\write\gtoutfile{To: math@arxiv.org}
\immediate\write\gtoutfile{Subject: put OR rep NNNNN:ppppp}
\immediate\write\gtoutfile{--text follows this line--}
\immediate\write\gtoutfile{Proxy-for: \ifx\theasciiauthors\relax
\theauthors\else\theasciiauthors\fi\s<\ifx\theasciiemail\relax\theemail\else\theasciiemail\fi>}
\immediate\write\gtoutfile{\noexpand\\}
\immediate\write\gtoutfile{Authors: \ifx\theasciiauthors\relax
\theauthors\else\theasciiauthors\fi}
{\def\\{ }\immediate\write\gtoutfile{Title: \ifx\theasciititle\relax
\thetitle\else\theasciititle\fi}}
\immediate\write\gtoutfile{Subj-class: GT or SG, GR etc}
\immediate\write\gtoutfile{MSC-class: \theprimaryclass\ifx\thesecondaryclass\relax\else, \thesecondaryclass\fi}
\immediate\write\gtoutfile{Journal-ref: Geom. Topol. Monogr. \thevolumenumber\s
(\thevolumeyear) \startpage-\finishpage}
\immediate\write\gtoutfile{Comments: Published by Geometry and Topology Monographs at}
\immediate\write\gtoutfile{\s\s\s  http://www.maths.warwick.ac.uk/gt/GTMon\thevolumenumber/paper\thepapernumber.abs.html}
\immediate\write\gtoutfile{\noexpand\\}
\immediate\write\gtoutfile{}
\ifx\theasciiabstract\relax
\immediate\write\gtoutfile{\theabstract}\else
\immediate\write\gtoutfile{\theasciiabstract}\fi
\immediate\write\gtoutfile{}
\immediate\write\gtoutfile{\noexpand\\}
\immediate\write\gtoutfile{}
\immediate\closeout\gtoutfile}}  %%% end of definition of \makeheadfile

\def\maketitlepage{\maketitlep\makeheadfile}
\let\makeshorttitle\maketitlepage
\let\maketitle\maketitlepage

\volumenumber{4}
\volumename{Invariants of knots and 3-manifolds (Kyoto 2001)}
\volumeyear{2002}
\papernumber{16}  
\pagenumbers{245}{261} 
\received{19 December 2001}
\revised{1 August 2002}
\accepted{10 September 2002}
\published{13 October 2002}

\usepackage{amssymb}
\usepackage{amsmath}
\usepackage{verbatim}
\usepackage{epsf}

\begingroup\makeatletter\ifx\SetFigFont\undefined%
\gdef\SetFigFont#1#2#3#4#5{%
  \reset@font\fontsize{#1}{#2pt}%
  \fontfamily{#3}\fontseries{#4}\fontshape{#5}%
  \selectfont}%
\fi\endgroup

\theoremstyle{plain}
\newtheorem*{theorem}{Theorem}
\newtheorem{prob}{Problem}

\theoremstyle{definition}

\newcommand{\bp}{\begin{prob}}
\newcommand{\ep}{\end{prob}}
\newcommand{\rks}{\rk{Remarks}}

\DeclareMathOperator{\Hom}{Hom}
\DeclareMathOperator{\Inv}{Inv}

\DeclareMathOperator{\vol}{vol}

\DeclareMathOperator{\td}{td}

\newcommand{\C}{\mathbb C}
\newcommand{\R}{\mathbb R}
\newcommand{\Z}{\mathbb Z}
\newcommand{\Q}{\mathbb Q}
\newcommand{\N}{\mathbb N}

\newcommand{\Rthree}{{\mathbb R}^3}

\newcommand{\ra}{\rightarrow}

\newcommand{\lan}{\langle}
\newcommand{\ran}{\rangle}

\newcommand{\LL}{\mathcal L}

\newcommand{\lie}[1]{\mathfrak{#1}}
\newcommand{\colie}[1]{\mathfrak{#1}^*}

\newcommand{\sixj}{\left\{ \begin{matrix} a&b&c\\ d&e&f
\end{matrix} \right\} }
\newcommand{\ksixj}{\left\{ \begin{matrix} ka&kb&kc\\ kd&ke&kf
\end{matrix} \right\} }
\newcommand{\kksixj}{\left\{ \begin{matrix} k\alpha&k\beta&k\gamma\\
k\delta&k\epsilon&k\zeta\end{matrix} \right\} }

\newcommand{\M}{\mathcal M}

\newenvironment{Relax}{\relax}{\relax}

\begin{document}

%%%%%%%%%%%%%%%%%%%%%%%%%%%%%%%%%%%%%%%%%%%%%%%%%%%%%%%%%%%%%%%%%%%%%
%%%%%%%%%%%%%%%%%%%%%%%%%%%%%%%%%%%%%%%%%%%%%%%%%%%%%%%%%%%%%%%%%%%%%

\title{Asymptotics and $6j$-symbols} 
\author{Justin Roberts}
\address{Department of Mathematics, UC San Diego\\9500 Gilman Drive,
La Jolla, CA 92093, USA}

\email{justin@math.ucsd.edu}
\primaryclass{22E99}\secondaryclass{81R05, 51M20}
\keywords{$6j$-symbol, asymptotics, quantization}

\begin{abstract}
Recent interest in the Kashaev-Murakami-Murakami hyperbolic volume
conjecture has made it seem important to be able to understand the
asymptotic behaviour of certain special functions arising from
representation theory --- for example, of the quantum $6j$-symbols for
$SU(2)$. In 1998 I worked out the asymptotic behaviour of the {\em
classical} $6j$-symbols, proving a formula involving the geometry of a
Euclidean tetrahedron which was conjectured by Ponzano and Regge in
1968. In this note I will try to explain the methods and philosophy
behind this calculation, and speculate on how similar techniques might
be useful in studying the quantum case.
\end{abstract}
\asciiabstract{
Recent interest in the Kashaev-Murakami-Murakami hyperbolic volume
conjecture has made it seem important to be able to understand the
asymptotic behaviour of certain special functions arising from
representation theory -- for example, of the quantum 6j-symbols for
$SU(2)$. In 1998 I worked out the asymptotic behaviour of the
classical 6j-symbols, proving a formula involving the geometry of a
Euclidean tetrahedron which was conjectured by Ponzano and Regge in
1968. In this note I will try to explain the methods and philosophy
behind this calculation, and speculate on how similar techniques might
be useful in studying the quantum case.}

\maketitle

\begin{Relax}\end{Relax}
%%%%%%%%%%%%%%%%%%%%%%%%%%%%%%%%%%%%%%%%%%%%%%%%%%%%%%%%%%%%%%%%%%%%%
%%%%%%%%%%%%%%%%%%%%%%%%%%%%%%%%%%%%%%%%%%%%%%%%%%%%%%%%%%%%%%%%%%%%%

\section{Introduction}

%%%%%%%%%%%%%%%%%%%%%%%%%%%%%%%%%%%%%%%%%%%%%%%%%%%%%%%%%%%%%%%%%%%%%
%%%%%%%%%%%%%%%%%%%%%%%%%%%%%%%%%%%%%%%%%%%%%%%%%%%%%%%%%%%%%%%%%%%%%

The Kashaev-Murakami-Murakami hyperbolic volume conjecture \cite{O,
Kashaev, Mur-Mur} is a conjecture about the asymptotic behaviour of a
certain sequence of ``coloured Jones polynomial'' knot invariants
$J_N(K)$, indexed by natural numbers $N$. In its simplest form, it
states that if the knot $K$ is hyperbolic, then the invariants grow
exponentially, with growth rate equal to the hyperbolic volume divided
by $2\pi$. We do not yet have any conceptual explanation of why this
conjecture might be true, and this seems a serious impediment to
attempts to prove it, despite the progress of Thurston \cite{T},
Yokota \cite{Y}, etc.

Attempts to prove and generalise this conjecture have led to renewed
interest in the asymptotics of the {\em quantum $6j$-symbols} for
$SU(2)$ and of the closely-related {\em Witten-Reshetikhin-Turaev} and
{\em Turaev-Viro} invariants of $3$-manifolds. The hope is that each
of these will display asymptotic behaviour governed by geometry in an
interesting and useful way.

What I want to describe in this note is a {\em philosophy}, a method
by which results of this form might be proved. In 1998 I proved an
essentially similar statement relating the asymptotic behaviour of the
{\em classical $6j$-symbols} to the geometry of Euclidean tetrahedra
\cite{R}. This theorem had been conjectured in 1968 by the physicists
Ponzano and Regge \cite{PR}; while well-known to and much used by
physicists, it had remained unproven and largely unexplained.

The method is {\em geometric quantization}: the idea is that if we
want to understand the asymptotic behaviour of some kind of
representation-theoretic quantity, then we should first write is as an
integral over some {\em geometrically meaningful} space, and use the
method of stationary phase to evaluate it in terms of local
contributions from (geometrically meaningful) critical points.

The plan of the paper is as follows. I start by defining $6j$-symbols
algebraically and describing some formulae for them. I then explain
their heuristic physical interpretation, and how this enabled Wigner
to give a rough asymptotic formula for them. I describe the general
method of geometric quantization in representation theory, with
special reference to the classical $6j$-symbol example. Finally I
explain how, at least in principle, one should be able to adapt these
techniques to deal with the quantum $6j$-symbol. I have tried to
complement rather than overlap the paper \cite{R} as much as
possible.

%%%%%%%%%%%%%%%%%%%%%%%%%%%%%%%%%%%%%%%%%%%%%%%%%%%%%%%%%%%%%%%%%%%%%
%%%%%%%%%%%%%%%%%%%%%%%%%%%%%%%%%%%%%%%%%%%%%%%%%%%%%%%%%%%%%%%%%%%%%

\section{The algebra of $6j$-symbols}

%%%%%%%%%%%%%%%%%%%%%%%%%%%%%%%%%%%%%%%%%%%%%%%%%%%%%%%%%%%%%%%%%%%%%
%%%%%%%%%%%%%%%%%%%%%%%%%%%%%%%%%%%%%%%%%%%%%%%%%%%%%%%%%%%%%%%%%%%%%

%\subsection{General definition of $6j$-symbols}

Suppose we have a category, such as the category of representations of
a compact group, possessing reasonable notions of tensor product,
duality, and decomposition into irreducibles. Let $I$ be a set
indexing the irreps, and let $V_a$ denote the irrep corresponding to
$a \in I$. Then there is an isomorphism
\[ V_a \otimes V_b \cong \bigoplus_{c \in I} V_c \otimes  \Hom(V_c, V_a
\otimes V_b)\] describing the decomposition of a tensor product of two
irreps. The space $\Hom(V_c, V_a \otimes V_b)$
might also be written as $\Inv(V_c^* \otimes V_a \otimes V_b)$, a
space of trilinear invariants or {\em multiplicity space}. With this
rule we can express arbitrary spaces of invariants in terms of
trilinear ones. Two obvious ways of decomposing a space of $4$-linear
invariants are 
\[ \Inv(V_a \otimes V_b \otimes V_c \otimes V_d) \cong \bigoplus_{e
\in I}
\Inv(V_a \otimes V_b \otimes V_e) \otimes \Inv(V_e^* \otimes V_c
\otimes V_d)\]
\[ \Inv(V_a \otimes V_b \otimes V_c \otimes V_d) \cong \bigoplus_{f
\in I}
\Inv(V_a \otimes V_c \otimes V_f) \otimes \Inv(V_f^* \otimes V_b
\otimes V_d)\] and the $6j$-symbol \[\sixj\] is defined to be the
part of the resulting ``change-of-basis'' operator mapping

$\Inv(V_a \otimes V_b \otimes V_e) \otimes \Inv(V_e^* \otimes V_c
\otimes V_d) \ra$\nl\hbox{}\hfill$ \Inv(V_a \otimes V_c \otimes V_f) 
\otimes
\Inv(V_f^* \otimes V_b \otimes V_d).$

The two most important properties of $6j$-symbols are their {\em tetrahedral
symmetry} and the {\em Elliott-Biedenharn} or {\em pentagon} identity.
The tetrahedral symmetry is a kind of equivariance property under
permutation of the six labels, summarised by associating it with a
labelled Mercedes badge: 
\[ \vcenter{\hbox{\mbox{\begin{picture}(0,0)%
\epsfbox{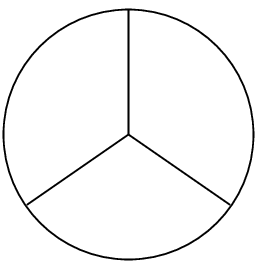}%
\end{picture}%
\setlength{\unitlength}{1973sp}%
\begin{picture}(2430,2436)(5086,-2175)
\put(5326,-586){\makebox(0,0)[lb]{\smash{\SetFigFont{10}{12.0}{\familydefault}{\mddefault}{\updefault}$f$}}}
\put(5851,-1561){\makebox(0,0)[lb]{\smash{\SetFigFont{10}{12.0}{\familydefault}{\mddefault}{\updefault}$b$}}}
\put(6226,-2086){\makebox(0,0)[lb]{\smash{\SetFigFont{10}{12.0}{\familydefault}{\mddefault}{\updefault}$d$}}}
\put(6676,-1561){\makebox(0,0)[lb]{\smash{\SetFigFont{10}{12.0}{\familydefault}{\mddefault}{\updefault}$c$}}}
\put(6001,-436){\makebox(0,0)[lb]{\smash{\SetFigFont{10}{12.0}{\familydefault}{\mddefault}{\updefault}$a$}}}
\put(7126,-586){\makebox(0,0)[lb]{\smash{\SetFigFont{10}{12.0}{\familydefault}{\mddefault}{\updefault}$e$}}}
\end{picture}
}}}\]
The Elliott-Biedenharn identity expresses the
fact that the composition of five successive change-of-basis operators
inside a space of $5$-linear invariants is the identity. For further
details on all of this see Carter, Flath and Saito \cite{CFS}.

%%%%%%%%%%%%%%%%%%%%%%%%%%%%%%%%%%%%%%%%%%%%%%%%%%%%%%%%%%%%%%%%%%%%%

%\subsection{The classical $6j$-symbols for $SU(2)$}

For the group $SU(2)$, things can be made much more concrete. Let $V$
denote the fundamental representation on $\C^2$, so that the
irreducible representations of $G$ are the symmetric powers $V_a =
S^aV$ ($a=0,1,2, \ldots$) with dimensions $a+1$. They are all {\em
self-dual}: $V_a \cong V_a^*$.

The spaces of trilinear invariants $\Inv(V_a \otimes V_b \otimes
V_c)$ are either one-dimensional or zero-dimensional,
according to whether $a,b,c,$ satisfy the following condition
``$(\Delta)$'' or not:
\[ a\leq b+c\qquad  b\leq c+a \qquad c\leq a+b \qquad a+b+c\ \hbox{is
even}.\] The triangle inequality here is the simplest example of the
``geometry governs algebra'' phenomenon with which we are concerned,
and it will be fully explained later.

Because these non-zero multiplicity spaces are one-dimensional, the
$6j$-symbols for $SU(2)$ are maps between one-dimensional vector
spaces, so by means of a suitable normalisation convention we can
think of them as numbers (in fact, they turn out to be {\em real}
numbers) rather than operators. By defining its value to be zero if
any of the triples don't satisfy $(\Delta)$, we can think of the
$6j$-symbol for $SU(2)$ as simply a real-valued function of six
natural numbers, which is invariant under the group $S_4$ of
symmetries of a tetrahedron.

There are various formulae for the $6j$-symbols. The {\em spin
network} method gives a straightforward but impractical combinatorial
formula. There is a one-variable summation of ratios of factorials,
which is the most efficient. There is also a generating function
approach. See Westbury \cite{Wes}, or \cite{CFS}.

%%%%%%%%%%%%%%%%%%%%%%%%%%%%%%%%%%%%%%%%%%%%%%%%%%%%%%%%%%%%%%%%%%%%%

%\subsection{Quantum $6j$-symbols}

The representation theory of the corresponding quantum group (Hopf
algebra) $U_q(\lie{sl}(2))$ has all the properties needed for
definition of $6j$-symbols, and has the same indexing of irreducibles,
resulting in the $\Q(q)$-valued {\em quantum $6j$-symbols} defined by
Kirillov and Reshetikhin \cite{KR}, which specialise at $q=1$ to the
classical ones. At a root of unity $q=e^{2\pi i /r}$ the
representation category may be quotiented to obtain one with finitely
many irreducibles, indexed by $0, 1, \ldots, r-2$. Turaev and Viro
\cite{TV} used the (real-valued) quantum $6j$-symbols associated to
this category to make an invariant of $3$-manifolds, and this is the
main reason for topologists to be interested in $6j$-symbols.

%%%%%%%%%%%%%%%%%%%%%%%%%%%%%%%%%%%%%%%%%%%%%%%%%%%%%%%%%%%%%%%%%%%%%

%\subsection{Application: Turaev-Viro state sum}

Let $T$ be a triangulation $T$ of a closed $3$-manifold $M$.  Define a
{\em state} $s$ to be an assignment of numbers in the range $0, 1,
\ldots, r-2$ to the edges of $T$. Given a state $s$, we can associate
to an edge $e$ labelled $s(e)$ the {\em quantum dimension} $d(s(e))$
of the associated irrep, and to each tetrahedron $t$ the {\em quantum
$6j$-symbol} $\tau(s(t))$ corresponding to the labels on its edges.
The real-valued state-sum
\[ Z(T) = \sum_s \prod_e d(s(e)) \prod_t \tau(s(t))\]
is invariant under the $2-3$ Pachner move, because of the
Elliott-Biedenharn identity. A minor renormalisation brings invariance
under the $1-4$ move too and so we obtain an invariant of $M$, the {\em
Turaev-Viro invariant} at $q =e^{2\pi i /r}$. 

The TV invariant turns out to be the square of the modulus of the
surgery-based {\em Witten-Reshetikhin-Turaev} invariant of $M$, which
therefore contains more information. But because it is computable in
terms of intrinsic structure (a triangulation), it should be easier to
relate to the geometry of $M$.

%%%%%%%%%%%%%%%%%%%%%%%%%%%%%%%%%%%%%%%%%%%%%%%%%%%%%%%%%%%%%%%%%%%%%
%%%%%%%%%%%%%%%%%%%%%%%%%%%%%%%%%%%%%%%%%%%%%%%%%%%%%%%%%%%%%%%%%%%%%

\section{The physics of classical $6j$-symbols}

%%%%%%%%%%%%%%%%%%%%%%%%%%%%%%%%%%%%%%%%%%%%%%%%%%%%%%%%%%%%%%%%%%%%%
%%%%%%%%%%%%%%%%%%%%%%%%%%%%%%%%%%%%%%%%%%%%%%%%%%%%%%%%%%%%%%%%%%%%%

%\subsection{Physical meaning}

To a physicist, the representation $V_a$ of $SU(2)$ is the space of
states of a quantum particle with spin $\frac12a$.  A {\em composite}
system of (for example) four particles with spins $\frac12a, \frac12b,
\frac12c, \frac12d$ is described by the {\em tensor product} of state
spaces. On this space there is a {\em total spin} operator (the
Casimir for the diagonal $SU(2)$, in fact) whose eigenspaces are the
irreducible summands; thus, for example, the invariant space $\Inv(V_a
\otimes V_b \otimes V_c \otimes V_d)$ is the subspace of states of the
system in which the total spin is zero.

The action of $SU(2)$ on the {\em first two} factors commutes with the
total spin operator, is and {\em its} Casimir gives the decomposition
\[ \Inv(V_a \otimes V_b \otimes V_c \otimes V_d) \cong \bigoplus_e
\Inv(V_a \otimes V_b \otimes V_e) \otimes \Inv(V_e \otimes V_c \otimes
V_d),\] into states in which the total spin of the first two (and
therefore also last two) particles is $\frac12e$.

The similar Casimir for the first and third particles does {\em
not} commute with this one and so gives a different eigenspace
decomposition. Standard quantum mechanics principles imply that the
square of the relevant matrix element
\[\sixj^2,\] 
is the {\em probability}, starting with the system in the state where
the first two particles have total spin $\frac12e$, that measuring the
total spin of the first and third combined gives $\frac12f$.

%%%%%%%%%%%%%%%%%%%%%%%%%%%%%%%%%%%%%%%%%%%%%%%%%%%%%%%%%%%%%%%%%%%%%

%\subsection{Wigner's semi-classical approximation}

The possible states of a classical particle with angular momentum of
magnitude $j$ are the vectors in $\R^3$ of length $j$. A {\em random}
such particle therefore has a state represented by a
rotationally-symmetric probability distribution on $\Rthree$ supported
on a sphere of radius $j$. For a {\em quantum} particle of a given
spin $j$, one can imagine the space of states as a space of certain
complex-valued wave-functions on $\Rthree$, whose pointwise norms give
(in general, rather spread-out) probability distributions for the
value of a hypothetical angular momentum vector. The {\em
semi-classical limit} requires that quantum particles with very large
spin should have distributions very close to those of the
corresponding classical particles, becoming more and
more localised near the appropriate sphere in $\Rthree$.

Wigner \cite{Wig} gave an asymptotic formula for the $6j$-symbols by adopting
this point of view. The {\em classical} version of the experiment
described above, whose output is the square of the $6j$-symbol, is as
follows. Suppose one has four random vectors of lengths $\frac12a,
\frac12b, \frac12c, \frac12d$ which form a closed quadrilateral; given
that one diagonal is $\frac12e$, what is the probability (density) that
the other is $\frac12f$? This analysis yielded the formula
\[ \sixj ^2 \approx \frac{1}{3 \pi V}, \]
with $V$ the volume of the Euclidean tetrahedron with edge-lengths $a,
b, \ldots, f$, supposing it exists. It should be
taken as a local root-mean-square average over the rapidly oscillatory
behaviour of the $6j$-symbol.

%%%%%%%%%%%%%%%%%%%%%%%%%%%%%%%%%%%%%%%%%%%%%%%%%%%%%%%%%%%%%%%%%%%%%

%\subsection{Ponzano-Regge state-sum}

There is a classical version of the Turaev-Viro state-sum, using edges
labelled by arbitrary irreps of $SU(2)$, which was written down by the
physicists Ponzano and Regge \cite{PR} in 1968. Their version is an
{\em infinite} state-sum which turns out to diverge for closed
$3$-manifolds; the Turaev-Viro invariant can be viewed as a
successful ``regularisation'' of their sum.

Their state-sum is a lattice model of {\em Euclidean quantum gravity},
which involves a path integral over the space of all Riemannian
metrics on a $3$-manifold. The states are interpreted as {\em
piecewise-Euclidean metrics} on $T$, made by gluing Euclidean
tetrahedra along faces, and from the asymptotic formula for
$6j$-symbols (below) one sees that the ``integrand'' measures the
curvature of the metric at the edges of $T$. Stationary points of the
``integral'' (classical solutions) should be metrics in which the
dihedral angles of the Euclidean simplexes glued around every edge sum
to $2\pi$, or at least to multiples of $2\pi$. (The resulting
ramification does seem to cause some problems in this model.)

Remarkably, the Turaev-Viro invariant with $q=e^{2\pi i/r}$ can be
interpreted in this context as a lattice model of quantum gravity {\em
with a positive cosmological constant}. Its stationary points should
correspond to metrics with {\em constant positive curvature}, and so
we should expect that the asymptotic behaviour of the TV invariant
(and likewise of the quantum $6j$-symbols themselves) as $r \ra
\infty$ will reflect this. For further details see the survey by Regge
and Williams \cite{RW}. Additional insight into the state-sum can be
obtained from Witten's paper \cite{Wit} or the work of Dijkgraaf and
Witten \cite{DW}.

%%%%%%%%%%%%%%%%%%%%%%%%%%%%%%%%%%%%%%%%%%%%%%%%%%%%%%%%%%%%%%%%%%%%%

%\subsection{Asymptotics}

We can associate to the six labels $a, b, \ldots f$ a metric
tetrahedron $\tau$ with these as side lengths.  The conditions
$(\Delta)$ guarantee that the individual faces may be realised
in Euclidean $2$-space, but as a whole the tetrahedron has an
isometric embedding into {\em Euclidean} or {\em Minkowskian}
$3$-space according to the sign of the {\em
Cayley determinant}, a cubic polynomial in the squares of the
edge-lengths. If $\tau$ is Euclidean, let $\theta_a, \theta_b, \ldots,
\theta_f$ be its corresponding exterior dihedral angles and $V$ its
volume.

\begin{theorem}{\rm\cite{R}}\qua
As $k \ra \infty$ {\rm (}for $k \in \Z${\rm )} there is an asymptotic formula
\[\ksixj \sim 
\begin{cases}{\displaystyle\sqrt{\frac{2}{3\pi Vk^3}}\cos{\left\{ \sum (ka+1)
\frac{\theta_a}{2} + \frac{\pi}{4}\right\}}} &\hbox{{if $\tau$ is
Euclidean,}} \\ \hbox{{exponentially decaying}}&\hbox{{if $\tau$ is
Minkowskian.}}\end{cases} \] {\rm (}The sum is over the
six edges of the tetrahedron.{\rm )}
\end{theorem}

To have a hope of proving this one needs to start from the {\em right
formula} for the $6j$-symbol. An approach very much in the spirit of
Wigner's is explained next.

%%%%%%%%%%%%%%%%%%%%%%%%%%%%%%%%%%%%%%%%%%%%%%%%%%%%%%%%%%%%%%%%%%%%%
%%%%%%%%%%%%%%%%%%%%%%%%%%%%%%%%%%%%%%%%%%%%%%%%%%%%%%%%%%%%%%%%%%%%%

\section{The geometry of classical $6j$-symbols}

%%%%%%%%%%%%%%%%%%%%%%%%%%%%%%%%%%%%%%%%%%%%%%%%%%%%%%%%%%%%%%%%%%%%%
%%%%%%%%%%%%%%%%%%%%%%%%%%%%%%%%%%%%%%%%%%%%%%%%%%%%%%%%%%%%%%%%%%%%%

%\subsection{Kirillov's method of orbits}

{\em Geometric quantization} is a collection of procedures for turning
symplectic manifolds (classical phase spaces) into Hilbert spaces
(quantum state spaces).  We will here
consider {\em K\"ahler quantization} only.

%\subsection{K\"ahler quantization}

If $M$ is a symplectic manifold with an {\em integral} symplectic form
(one that evaluates to an integer on all classes in $H_2(M;\Z)$) then
it is possible to find a smooth line bundle $\LL$ on $M$ with a
connection whose curvature form is $(-2\pi i)^{-1}\omega$. The
quantization $Q(M)$ is then a subspace of the space of sections of
$\LL$, specified by a choice of {\em polarisation} of $M$.

If $M$ is {\em K\"ahler} (complex in a way compatible with the
symplectic form) then there is a standard way to polarise it: the
bundle $\LL$ can be taken to be {\em holomorphic}, and the
relevant subspace $Q(M)$ is its space of {\em holomorphic sections}.

Such a bundle can also be given a smooth {\em hermitian metric} $\lan
-, -\ran$ compatible with its connection. When $M$ is compact, the
space $Q(M)$ will be finite-dimensional, and we can define an obvious
Hilbert space inner product of two sections by the integral formula
\[ ( s_1, s_2 ) = \int_M \lan s_1, s_2 \ran \frac{\omega^n}{n!}. \]
The {\em dimension} of $Q(M)$ can be computed cohomologically via the
{\em Riemann-Roch formula}: at least, the {\em Euler characteristic}
$\chi(\LL)$ of the set of sheaf cohomology groups $H^*(M;\LL)$ is
given by
\[ \int_M e^{c_1(M)} \td(TM),\]
and in many cases one can prove a {\em vanishing theorem} showing that
the space of holomorphic sections $H^0(M;\LL)$ is the only non-trivial
space, and thereby obtain a direct formula for its dimension. 

Note that we can rescale the symplectic form by a factor of $k \in
\N$, replacing $\LL$ by $\LL^{\otimes k}$, and repeat the
construction. Examining the behaviour as $k \ra \infty$ corresponds to
examining the behaviour of the quantum system as $\hbar = 1/k$ tends
to zero; this is the {\em semi-classical limit}. If $M$ has dimension
$2n$ then the formula for $\chi(\LL^{\otimes k})$ is a polynomial in
$k$ with leading term $k^n \vol(M)$, where the volume is measured with
respect to the symplectic measure $\omega^n/n!$. This phenomenon is
the simplest possible manifestation of the kind of geometric
asymptotic behaviour we are studying.

Note also that if there is an equivariant action of a compact group
$G$ on $\LL \ra M$ which preserves the K\"ahler structure and
hermitian form then it acts on the sections of $\LL$, giving a {\em
unitary} representation of $G$. In this case there is an {\em
equivariant index formula} giving the character of $H^*(M;\LL)$ in
cohomological terms, and also a {\em fixed-point formula} for the
character which may be regarded as a kind of {\em exact}
semi-classical approximation.

%\subsection{Borel-Weil-Bott theorem}

If $G$ is a Lie group with Lie algebra $\lie{g}$ then its {\em
coadjoint representation} $\lie{g}^*$ decomposes as a union of
symplectic {\em coadjoint orbits} under the action of $G$. Kirillov's
{\em orbit principle} \cite{K2} is that quantization induces a
correspondence between the irreducible unitary representations of $G$
and certain of the coadjoint orbits, though the association does
depend on the method of quantization used. For a {\em compact} group
$G$, the coadjoint orbits are $G$-invariant K\"ahler manifolds, and the
{\em integral} ones are parametrised by (in fact, are the orbits
through) the weights in the positive Weyl chamber. K\"ahler
quantization turns the orbit through the weight $\lambda$ into the
irrep with highest weight $\lambda$. This is (part of) the {\em
Borel-Weil-Bott theorem}, which is described more
algebro-geometrically in Segal \cite{CSM} or Fulton and Harris
\cite{FH}.

%\subsection{Invariant theory}

The correspondence between K\"ahler manifolds and representations is
very helpful in understanding invariant theory for Lie groups. There
are three essential ideas: first, the above association between irreps
and integral coadjoint orbits; second, that tensor products of
representations correspond to products of K\"ahler manifolds; third,
that taking the space of $G$-invariants of a representation
corresponds to taking the {\em K\"ahler quotient} of a manifold. 

The $G$-actions we are dealing with are {\em Hamiltonian}, meaning
that the vector fields defining the infinitesimal action of $G$ are
symplectic gradients and that we can define an equivariant {\em moment
map} $\mu: M \ra \colie{g}$ collecting them all up according to the
formula
\[ d\mu(\xi) = \iota_{X_\xi}\omega \quad (= \omega(X_\xi, -)), \]
where $\xi \in \lie{g}$ and $X_\xi$ is the corresponding vector
field. The K\"ahler quotient is then defined as
$M_G=\mu^{-1}(0)/G$. The theorem of Guillemin and Sternberg \cite{GS}
is that $Q(M_G) = \Inv(Q(M))$. Note that for a coadjoint orbit the
moment map turns out to be simply the inclusion map $M \subseteq
\colie{g}$, and for a product of manifolds, the moment map is the sum of
the individual ones. 

%\subsection{The example of $SU(2)$}

For $SU(2)$ the coadjoint space is Euclidean $\Rthree$, with the group
acting by $SO(3)$ rotations. All coadjoint orbits other than the
origin are spheres, and the integral ones are those $S^2_a$ with
integral radius $a$. K\"ahler quantization entails thinking of $S^2_a$
as the Riemann sphere, equipped with the $a$th tensor power of the
hyperplane bundle; the space of holomorphic sections is the irrep
$V_a$, and Riemann-Roch gives its dimension (correctly!) as $a+1$.

To compute the space $\Inv(V_a \otimes V_b \otimes V_c)$, we first
form the product $M$ of the three spheres of radii $a,b,c$. Its moment
map is just the sum of the three inclusion maps into $\Rthree$, so
that $\mu^{-1}(0)$ is the space of closed triangles of vectors of
lengths $a,b,c$. Now $M_G$ is the space of such things up to overall
rotation: it is either a point or empty, and its quantization
$Q(M_G)=\Inv(V_a \otimes V_b \otimes V_c)$ is either $\C$ or zero,
according to the triangle inequalities, whose role in $SU(2)$
representation theory is now apparent. (The additional parity
condition can only be seen by considering the lift of the $SU(2)$
action to the line bundle $\LL$.) Higher ``polygon spaces'' arise
similarly: for example, $\Inv(V_a \otimes V_b \otimes V_c \otimes V_d)$
is the quantization of the moduli space of shapes of quadrilaterals of
sides $a,b,c,d$ in $\R^3$.

%\subsection{Localisation of invariant sections}

A fundamental ingredient of Guillemin and Sternberg's proof that
quantization commutes with reduction is the fact that a $G$-invariant
section $s$ of the equivariant bundle $\LL \ra M$ has maximal
pointwise norm on the set $\mu^{-1}(0)$. In fact, the norm of $s$
decays in a Gaussian exponential fashion in the transverse directions 
(and will in fact reach the value zero on the {\em unstable points} of
$M$).

The $k$th power $s^k$ is an invariant section of $\LL^{\otimes k}$
whose norm decays faster; we can imagine in the limiting case $k \ra
\infty$ that such a section becomes {\em localised} to a
delta-function-like distribution supported on $\mu^{-1}(0)$. {\em
Pairings} of such sections will become localised to the {\em
intersections} of these support manifolds, and this is the main idea
of the proof of the asymptotic formula for the $6j$-symbol. 

In \cite{R} it is written as a pairing between two $12$-linear
invariants, and thus as an integral over the symplectic quotient of
the product of twelve spheres, whose radii depend on the six
labels. The intersection locus amounts either to two points
corresponding to mirror-image Euclidean tetrahedra, if they exist, or
is empty. In the first case one gets a sum of two local contributions,
each a Gaussian integral, and after rather messy calculations the
formula emerges; exponential decay is automatic in the second case.

%%%%%%%%%%%%%%%%%%%%%%%%%%%%%%%%%%%%%%%%%%%%%%%%%%%%%%%%%%%%%%%%%%%%%
%%%%%%%%%%%%%%%%%%%%%%%%%%%%%%%%%%%%%%%%%%%%%%%%%%%%%%%%%%%%%%%%%%%%%

\section{The geometry of quantum $6j$-symbols}

%%%%%%%%%%%%%%%%%%%%%%%%%%%%%%%%%%%%%%%%%%%%%%%%%%%%%%%%%%%%%%%%%%%%%
%%%%%%%%%%%%%%%%%%%%%%%%%%%%%%%%%%%%%%%%%%%%%%%%%%%%%%%%%%%%%%%%%%%%%

%\subsection{The ``curved version'' of the method of orbits}

Quantum $6j$-symbols, evaluated at a root of unity, come from a
category which might be considered as the category of representations
of a quantum group at this root of unity, or of a loop group at a
corresponding level. The loop group picture leads to a beautiful and
conceptually very valuable analogue of the geometric framework
described above. In principle it also allows an analogous calculation
of the asymptotic behaviour, though in practice this seems quite
difficult.

For a compact group $G$ we saw above the association between irreps
and integral coadjoint orbits. Let us now consider the analogous
correspondence for the category of positive energy representations of
its loop group $LG$ at level $k$. (For the actual construction of the
representations see Pressley and Segal \cite{PS}.)

Instead of coadjoint orbits we use {\em conjugacy classes} in $G$
itself. Notice that the foliation of $\lie{g}$ by adjoint orbits is a
{\em linearisation} of the foliation of $G$ by conjugacy classes at
the identity, so that the quantum orbit structure is a sort of {\em
curved counterpart} of the classical case. The conjugacy classes
correspond under the exponential map to points of a {\em Weyl alcove},
a truncation of a Weyl chamber inside a Cartan subalgebra. The
``integral'' conjugacy classes giving the irreps at level $k$ are
those obtained by exponentiating $k^{-1}$ times the elements of the
weight lattice lying in a $k$-fold dilation of this alcove.

The definition of the {\em fusion} tensor product of such irreps is
subtle. However, given integral conjugacy classes $C_1, C_2, \ldots,
C_n$ corresponding to level-$k$ irreps of $LG$, it is not hard to
describe a symplectic manifold $\M(C _1, C_2, \ldots, C_n)$ which will
correspond to the invariant part of their tensor product.

Let $\M$ be the moduli space of flat $G$-connections on an
$n$-punctured sphere. These are just representations, up to conjugacy,
of its fundamental group, which we take to have one generator for each
boundary circle and the relation that their product is $1$. This space
$\M$ is a Poisson manifold, and traces of the puncture holonomies give
Casimir functions on it. Their common level sets, the symplectic
leaves, are the spaces $\M(C _1, C_2, \ldots, C_n)$ comprising
representations with the generators mapping to given conjugacy
classes. 

In the case of $SU(2)$, the exponential map gives a bijective correspondence
between the unit interval in the Cartan subalgebra $\R$ and the
conjugacy classes. At level $k$, the allowable highest weights are
therefore $0, 1, \ldots, k$, corresponding to the conjugacy classes
$C_a$ with trace equal to $2 \cos (\pi a/k)$ for some $a=0, 1, \ldots,
k$, and to the positive energy irreps $V_a$ of $LSU(2)$ at level $k$.

The space of trilinear invariants $\Inv (V_a \otimes V_b \otimes V_c)$
(where $0\leq a,b,c \leq k$) corresponds to the space $\M(C_a, C_b,
C_c)$ of triples of matrices inside $C_a \times C_b \times C_c$ whose
product is $1$, considered up to conjugacy. This amounts to the space
of shapes of triangles of sides $\frac{a}{k}, \frac{b}{k},
\frac{c}{k}$ in {\em spherical} $3$-space, and is either a single
point or empty according to the {\em quantum triangle inequalities},
meaning the condition $(\Delta)$ together with the extra rule
$a+b+c\leq 2k$. This corresponds to the well-known fusion rule for
quantum $SU(2)$ at root of unity $q=e^{2 \pi i/(k+2)}$. Similarly,
spaces of quadrilinear invariants correspond to spaces of spherical
quadrilaterals with prescribed lengths, and so on. This is the real
justification for the use of the word ``curved'' above!

%\subsection{Quantization}

The explicit construction of a vector space (of tensor invariants)
from a symplectic manifold such as $\M(C_1, C_2, \ldots, C_n)$ is
achieved as before by using the K\"ahler quantization technique. The
technical difference here is that such spaces have {\em many} natural
complex structures, and so the procedure is more subtle. A choice of
complex structure on the underlying punctured sphere induces a complex
structure on $\M(C_1, C_2, \ldots, C_n)$ which can be used to
construct a holomorphic line bundle and a finite-dimensional space of
holomorphic sections, the {\em space of conformal blocks}. These
spaces depend smoothly on the chosen complex structure and in fact
form a bundle over the Teichm\"uller space of such structures with a
natural {\em projectively flat connection}, described by Hitchin
\cite{H}. The connection enables canonical and coherent
identications of all the different fibre spaces, at least up to
scalars.

As before, there is an index formula, the {\em Verlinde formula}, for
the dimensions of such spaces. We can consider the semi-classical
limit by sending the level $k$ to infinity but keeping the conjugacy
classes fixed, because the highest weight of the irrep corresponding
to a fixed conjugacy class scales with the level. The formula is then
a polynomial in $k$ with leading term given by the volume of $\M(C_1,
C_2, \ldots, C_n)$.  See Witten \cite{Wit} or Jeffrey and Weitsman
\cite{JW} for more details here, though perhaps the theory of {\em
quasi-Hamiltonian spaces} developed by Alexeev, Malkin and Meinrenken
\cite{AMM} will ultimately give the best framework.

%\subsection{Computing the quantum $6j$-symbol}

The asymptotic problem for the $SU(2)$ quantum $6j$-symbol is as
follows. Pick six rational numbers $\alpha, \beta, \gamma, \delta,
\epsilon, \zeta$ between $0$ and $1$. For a level $k$ such that the
the six products $a=\alpha k$ etc. are integers, we want to evaluate the
quantum $6j$-symbol
\[ \kksixj\] 
at $q=e^{2\pi i /(k+2)}$ and then look at the asymptotic expansion as
$k \ra \infty$. The guess is that this should have something to do
with the geometry of a {\em spherical} tetrahedron, since we have
everywhere replaced geometry of the original coadjoint $\R^3$ with the
group $SU(2)=S^3$.

We can express this quantum $6j$-symbol as a hermitian pairing between
a certain pair of vectors in the space $\Inv(V_a \otimes V_b \otimes
V_c \otimes V_d)$. This means working over a $2$-dimensional
symplectic manifold, the space of spherical quadrilaterals of given
edge-lengths.

Now the classical version of this manifold, the space of Euclidean
quadrilaterals, has well-known Hamiltonian circle actions
corresponding to the lengths of the diagonals of the
quadrilateral. This extends to the quantum, spherical case: the
lengths are in fact the traces of the holonomies around curves
separating the punctures into pairs, and generate Goldman's flows
\cite{JW}.

If such a circle action were to preserve the K\"ahler structure then
it would act on the quantization, thereby decomposing the space of
quadrilinear invariants into one-dimensional weight spaces. It would
be natural to assume that these would generate the different bases
mentioned in section 2 and the vectors we need to pair to compute the
$6j$-symbol.

In the classical case these flows do {\em not} preserve the natural
K\"ahler structure on the product of four spheres. To proceed one
would need additional machinery to show that the quantization is
independent of the K\"ahler polarisation; then one would recover the
action of the circles on the quantization and perhaps be able to carry
the idea through, obtaining an alternative to the proof in \cite{R}.

In the quantum case we have {\em moduli} of K\"ahler structures coming
from the choices of complex structure on the sphere with $4$
distinguished points. The moduli space is a Riemann sphere minus three
points; these correspond to ``stable curve'' degenerations and we may
add them in to compactify it.

At each singular point there is a {\em Verlinde decomposition} of the
space of conformal blocks into a sum of tensor products of
one-dimensional trilinear invariant spaces, and these spaces {\em are}
the eigenspaces of the Hamiltonian flow which preserves the degenerate
complex structure. So we ought to be able to specify geometrically the
two sections we need to pair. Unfortunately they live in the fibres
over different points in the moduli space, so we then need to parallel
transport one using the the holonomy of the projectively flat
connection before we can pair them easily. Dealing with this might be
difficult; it seems for example that even the unitarity of the holonomy is
still not established.

If we view the moduli space as $\C -\{0,1\}$ then we seek the holonomy
along the unit interval from $0$ to $1$. Now {\em asymptotically} the
connection we are examining becomes the {\em Knizhnik-Zamolodchikov
connection}, and this holonomy is nothing more than the {\em Drinfeld
associator}. (See Bakalov and Kirillov \cite{BK}, for example.) This
is the {\em geometric} explanation for the equivalence of the
$6j$-symbol and associator pointed out recently by Bar-Natan and
Thurston. Of course, one could try to compute a nice tetrahedrally
symmetric formula for the associator hoping that the asymptotic
formula for the $6j$-symbol would follow: this would be a completely
alternative approach to the asymptotic problem.

%%%%%%%%%%%%%%%%%%%%%%%%%%%%%%%%%%%%%%%%%%%%%%%%%%%%%%%%%%%%%%%%%%%%%
%%%%%%%%%%%%%%%%%%%%%%%%%%%%%%%%%%%%%%%%%%%%%%%%%%%%%%%%%%%%%%%%%%%%%

\section{Related problems}

%%%%%%%%%%%%%%%%%%%%%%%%%%%%%%%%%%%%%%%%%%%%%%%%%%%%%%%%%%%%%%%%%%%%%
%%%%%%%%%%%%%%%%%%%%%%%%%%%%%%%%%%%%%%%%%%%%%%%%%%%%%%%%%%%%%%%%%%%%%

\bp Compute the asymptotics of the quantum $6j$-symbol. \ep

\rks One programme for the computation was outlined above. Chris
Woodward has recently conjectured \cite{Woo} a precise formula and
checked it empirically. Suppose there is a spherical tetrahedron with
sides $l$ equal to $\pi$ times $\alpha, \beta, \gamma, \delta,
\epsilon, \zeta$, and associated dihedral angles $\theta_l$. Let $V$
be its volume, and let $G$ be the determinant of the spherical Gram
matrix, the symmetric $4\times 4$ matrix with ones on the diagonal
and the quantities $\cos (l)$ off the diagonal. Then he conjectures
that
\[ \kksixj_{q=e^{2\pi i /(k+2)}} \sim \sqrt{\frac{4 \pi^2}{k^3
\sqrt{G}}} 
\cos\left\{\sum (kl+1)
\frac{\theta_l}{2} - \frac{k}{\pi} V + \frac{\pi}{4}\right\}. \]

%%%%%%%%%%%%%%%%%%%%%%%%%%%%%%%%%%%%%%%%%%%%%%%%%%%%%%%%%%%%%%%%%%%%%

\bp Compute the asymptotics of the Turaev-Viro invariant of a closed
$3$-manifold. \ep

\rks Such a formula might result from an asymptotic formula for the
quantum $6j$-symbol, though this would not be straightforward.  The
hope is that the asymptotics might relate to the existence of
spherical geometries on a $3$-manifold, although technical issues
related to ramified gluings of spherical tetrahedra make this seem
likely to be a fairly weak connection.

%%%%%%%%%%%%%%%%%%%%%%%%%%%%%%%%%%%%%%%%%%%%%%%%%%%%%%%%%%%%%%%%%%%%%

\bp Prove the Minkowskian part of Ponzano and Regge's formula. \ep

\rks The methods of \cite{R} don't give a precise formula for the
exponentially decaying asymptotic regime which occurs when the
stationary points have become ``imaginary''. It is possible {\em
formally} to write down a ``Wick rotated'' integral over a product of
hyperboloids, instead of spheres, as a formula for the same classical
$6j$-symbol. This integral has well-defined stationary points
corresponding to Minkowskian tetrahedra, whose local contributions
seem correct, but the problem is that it does not converge! Some kind
of argument involving deformation of the contour of integration is
probably required.

%%%%%%%%%%%%%%%%%%%%%%%%%%%%%%%%%%%%%%%%%%%%%%%%%%%%%%%%%%%%%%%%%%%%%

\bp Try to compute asymptotic expansions of similar quantities. \ep

\rks Classical $6j$-symbols can be generalised to so-called
$3nj$-symbols, associated to arbitrary trivalent labelled graphs drawn
on a sphere. The asymptotic behaviour here will be governed by many
stationary points corresponding to the different isometric embeddings
of such a graph into $\R^3$, but it's not completely clear what the
expected contribution from each should be --- the volume, or something
more complicated.

Stefan Davids \cite{D} studied $6j$-symbols for the non-compact group
$SU(1,1)$. There are various different cases corresponding to unitary
irreps from the discrete and continuous series, and some surprising
relations between the discrete series symbols and geometry of
Minkowskian tetrahedra. 

One could study the Frenkel-Turaev elliptic and trigonometric
$6j$-symbols \cite{FT}. I have no idea what they might correspond to
geometrically.

The $6j$-symbols for higher rank groups are not simply scalar-valued
quantities, because the trilinear invariant spaces typically have
dimension bigger than one. This makes them trickier to handle and the
question of asymptotics less interesting. One could at least study
their norms as operators and expect a geometrical result. There are
possibly some nice special cases: Knutson and Tao \cite{KT} showed
that for $GL(N)$ the property of three irreps having multiplicity one
is stable under rescaling their highest weights by $k$, so one can
expect some scalar-valued $6j$-symbols with interesting asymptotics.

%%%%%%%%%%%%%%%%%%%%%%%%%%%%%%%%%%%%%%%%%%%%%%%%%%%%%%%%%%%%%%%%%%%%%

\bp The hyperbolic volume conjecture. \ep

\rks A basic approach to the conjecture is to try to give a formula
for Kashaev's coloured Jones polyomial invariant in terms of some
quantum dilogarithms associated to an ideal triangulation of the knot
complement, and then relate their asymptotics to geometry. In fact the
quantum dilogarithm, the basic ingredient in Kashaev's invariant, does
seem to behave as a kind of $6j$-symbol, satisfying a pentagon-type
identity and having an asymptotic relationship with volumes of ideal
hyperbolic tetrahedra. It would seem helpful to be able to interpret
it as arising from geometric quantization of some suitable space of
hyperbolic tetrahedra, with a view to gaining conceptual understanding
of the conjecture.

Jun Murakami and Yano \cite{JM} have applied Kashaev's non-rigorous
stationary phase methods to the Kirillov-Reshetikhin sum formula for
the quantum $6j$-symbol. The (false) result is exponential growth,
with growth rate given by the volume of the hyperbolic tetrahedron
with the appropriate dihedral angles. Similarly, Hitoshi Murakami
\cite{HM} has obtained ``fake'' exponential asymptotics for the
Turaev-Viro invariants of some hyperbolic $3$-manifolds.

These strange results are very interesting. In each case we start with
something which can be expressed as an $SU(2)$ path integral and in an
alternative way as a sum. Applying perturbation theory methods to the
path integral suggests the correct polynomial asymptotics.  But
``approximating'' the sum by a contour integral in the most obvious
way and applying stationary phase gives very different asymptotic
behaviour, seemingly reflecting a {\em complexification} of the
original path integral. There is a certain similarity to the
appearance of the ``imaginary'' Minkowskian critical points in problem
3. As in that case, the problem appears to be making sense of the
complexified integral in the first place. It is presumably this
quantity which we should be interested in as a genuine exponentially
growing invariant, and which we should try to learn to compute using
some kind of TQFT techniques.

\rk{Acknowledgements}\qua This note describes work carried out under an
EPSRC Advanced Fellowship, NSF Grant DMS-0103922 and JSPS fellowship
S-01037. 

%%%%%%%%%%%%%%%%%%%%%%%%%%%%%%%%%%%%%%%%%%%%%%%%%%%%%%%%%%%%%%%%%%%%%
%%%%%%%%%%%%%%%%%%%%%%%%%%%%%%%%%%%%%%%%%%%%%%%%%%%%%%%%%%%%%%%%%%%%%

\Addresses
\end{document}

%%%%%%%%%%%%%%%%%%%%%%%%%%%%%%%%%%%%%%%%%%%%%%%%%%%%%%%%%%%%%%%%%%%%%
%%%%%%%%%%%%%%%%%%%%%%%%%%%%%%%%%%%%%%%%%%%%%%%%%%%%%%%%%%%%%%%%%%%%%